\documentclass[12pt]{article}
\usepackage{amssymb,latexsym}
\usepackage{graphics}

\def\wrtext#1{\relax\ifmmode{\leavevmode\hbox{#1}}\else{#1}\fi}
\def\abs#1{\left|#1\right|}
\def\begeq{\begin{equation}}
\def\endeq{\end{equation}}
\def\Remark{\vskip 2mm \noindent {\em Remark}}

\def\neigh{neighborhood}

\def\Re{{\rm Re\,}}
\def\Im{{\rm Im\,}}

\newcommand{\eps}{\epsilon}
\def\part#1{\frac{\partial}{\partial #1}}

\newcommand{\real}{\mbox{\bf R}}
\newcommand{\comp}{\mbox{\bf C}}

\newcommand{\nat}{\mbox{\bf N}}

\renewcommand{\Re}{\mbox{\rm Re\,}}
\renewcommand{\Im}{\mbox{\rm Im\,}}

\renewcommand{\exp}{\mbox{\rm exp\,}}

\newtheorem{dref}{Definition}[section]

\newtheorem{theo}[dref]{Theorem}
\newtheorem{prop}[dref]{Proposition}

\newenvironment{proof}{\vspace{.3cm}\noindent{{\em Proof:}}}{\hfill$\Box$}
\begin{document}

\begin{center}
\Large {\bf Boundary spectral behaviour for semiclassical operators in
one dimension}
\end{center}

\bigskip
\noindent {\bf Michael Hitrik}

\medskip
\noindent Department of Mathematics, University of California, Los
Angeles, CA 90095-1555, USA \\
E-mail: hitrik@math.ucla.edu


\medskip
\noindent {\bf 1991/2000 AMS Subject Classification:} 34E05,
34L20, 47A12


\vspace*{1cm}
\noindent
{\bf Abstract}: For a class of non-selfadjoint semiclassical operators in dimension
one, we get a complete asymptotic description of all eigenvalues near a
critical value of the leading symbol of the operator on the boundary of the
pseudospectrum.

\section{Introduction and main result}\label{section0}
\setcounter{equation}{0}

In his work~\cite{Davies1}, E. B. Davies considers the
non-selfadjoint harmonic oscillator,
$$
P(h)=h^2 D_x^2+e^{i\alpha} x^2,\quad 0<\alpha<\pi/2.
$$
By an analytic continuation argument, it was established that the
spectrum of $P(h)$ consists of the eigenvalues
$\{e^{i\alpha/2}h(2k+1),\,\,k=0,1,2,\ldots\}$, while the
semiclassical pseudospectrum of $P(h)$, defined as the range of
the symbol of $P(h)$, $p(x,\xi)=\xi^2+e^{i\alpha} x^2$, on
$\real^2$, is equal to the sector $\{z\in \comp; 0\leq \arg z\leq
\alpha\}$. As $h\rightarrow 0$, the spectrum accumulates precisely
at the origin, which is the only point $z_0$ at the boundary of
the pseudospectrum, for which $dp=0$ on $p^{-1}(z_0)\cap \real^2$.
The angle at which the eigenvalues approach the origin is given
here by $\alpha/2$. Now associated with the quadratic form $p$,
there is a Hamilton map $F$, defined by
$$
F=\frac{1}{2}\left(\begin{array}{ccc}
  p''_{x \xi} & p''_{\xi \xi} \\
  -p''_{x x} & -p''_{x \xi}
  \end{array}\right)=\left(\begin{array}{ccc}
                       0 & 1 \\
                       -e^{i\alpha} & 0
                       \end{array}\right),
$$
with the eigenvalues $\pm \mu$, where $\mu=i e^{i\alpha/2}$. We
may therefore reformulate the result of~\cite{Davies1} by saying
that the angle at which the spectrum approaches the boundary point
$0$ is given by $\arg(\mu/i)$. In this paper we shall give a
straightforward extension of this result to what we think is a
natural class of non-selfadjoint $h$-pseudodifferential operators
in dimension one. We shall show that, in the semiclassical limit,
the direction of the spectrum near a critical value on the
boundary of the pseudospectrum is determined by a suitable
eigenvalue of the Hamilton map of the Hessian of the principal
symbol at a corresponding critical point.

Without any restriction on the dimension, pseudospectra of
non-selfadjoint semiclassical pseudodifferential operators have
been studied in the recent paper~\cite{DSZ}, which extended some of the
results of~\cite{Davies2},~\cite{Zw1},~\cite{Zw2}.
Following~\cite{Zw2} and~\cite{DSZ}, let us consider the
convection-diffusion operator with a quadratic potential,
$$
P(x,hD_x)=(hD_x)^2+i(hD_x)+x^2=e^{x/2h}\circ\left((hD_x)^2+x^2+\frac{1}{4}\right)\circ
e^{-x/2h}.
$$
Using the last representation, one can prove that the spectrum of
$P(x,hD_x)$ is given by $\{(2k+1)h+\frac{1}{4},\,k=0,1,2\ldots\}$,
while the range of the symbol of $P(x,hD_x)$ fills out the region
$\{z\in \comp; \Re z\geq (\Im z)^2\}$. In this case, the spectrum
lies strictly inside the pseudospectrum, and Theorem 3
in~\cite{DSZ} gives general conditions on a point $z_0$ at the
boundary of the pseudospectrum of an operator $P=p^w(x,hD_x)+{\cal
O}(h)$, to be away from the spectrum, as $h\rightarrow 0$. The
first one is the principal type condition, \begeq \label{1.0} dp\neq
0\,\,\,\wrtext{along}\,\,\,p^{-1}(z_0),
\endeq
and the second is the following dynamical hypothesis, \begeq
\label{1.0.1} \wrtext{For some}\,\lambda\in \comp,\,\wrtext{no
complete trajectory of}\, H_{{\rm Re} (\lambda p)}\,\wrtext{is
contained in}\,p^{-1}(z_0).
\endeq
Asssuming (\ref{1.0}) and (\ref{1.0.1}) together with a basic
analyticity assumption, it is established in~\cite{DSZ} that there
is a sufficiently small but fixed \neigh{} of $z_0$ in $\comp$
which does not intersect the spectrum of $P$. Back to the
one-dimensional case, the purpose of this note is to describe the
eigenvalue distribution near a boundary point $z_0$, in a simplest
situation when the condition (\ref{1.0}) fails to hold (so that
(\ref{1.0.1}) fails as well). In doing so, as in~\cite{DSZ}, we
shall work under the analyticity assumptions, which will be essential
for our methods.

\vskip 4mm \noindent
Let
$$
P=P^w(x,hD_x;h)
$$
be the $h$--Weyl quantization of a symbol $P(x,\xi;h)$, which is a
holomorphic function of $(x,\xi)$ in a tubular \neigh{} of
$\real^2\subset \comp^2$, with \begeq \label{1.1} P(x,\xi;h)={\cal
O}(1) m(\Re (x,\xi))
\endeq
there. Here we assume that $1\leq m\in C^{\infty}(\real^2)$ is an
admissible weight function in the sense that for some fixed
$C_0,N>0$ we have
$$
m(X)\leq C_0 \langle{X-Y}\rangle^{N} m(Y),\quad X,Y\in \real^2.
$$
We next assume that as $h\rightarrow 0$, $P(x,\xi;h)$ has an
expansion in the space of holomorphic functions satisfying the bound (\ref{1.1}),
$$
P(x,\xi;h)\sim\sum_{j=0}^{\infty} h^j p_j(x,\xi).
$$
For $h>0$ small enough, and when equipped with the domain $H(m)$,
the naturally defined Sobolev space associated to the weight $m$,
$P$ becomes a closed densely defined operator on $L^2(\real)$.

We shall denote the principal symbol of $P$ by $p:=p_0$, and make
an additional assumption of the ellipticity of $p$ near infinity,
\begeq \label{1.2}
 \abs{p(X)}\geq \frac{m(\Re X)}{C}, \quad
\abs{X}\geq C, \quad \abs{\Im X}\leq \frac{1}{C},\quad C>0.
\endeq
Throughout this paper,  we shall work under the assumption that
\begeq \label{1.3}
m(X)\rightarrow \infty,\quad \abs{X}\rightarrow
\infty,\,\,X\in \real^2.
\endeq
As before, associated to $P$, we introduce the semiclassical
pseudospectrum \begeq \label{1.3.5} \Sigma(p):=p(\real^2).
\endeq
This is a closed set, since $p$ is proper, and we shall assume
that $\Sigma(p)$ is not all of $\comp$. It follows then from
(\ref{1.3}) that if $z_0$ is in the complement of $\Sigma(p)$,
then $(P-z_0)^{-1}$ exists and is a compact operator on $L^2$. The
analytic Fredholm theory implies that the spectrum of
$P$ is discrete and consists of eigenvalues of finite
multiplicity.

Let $z_0\in \partial \Sigma(p)$ and assume that
$p^{-1}(z_0)\cap\real^2$ is a finite collection of points. It will
in fact be sufficient to treat the case when the pre-image of
$z_0$ is just a single point, say $(0,0)\in \real^2$, and in what
follows we shall work under this assumption. We introduce an
exterior cone condition, \begin{eqnarray} \label{1.3.6} & &
\wrtext{There exist}\,\,\eps_0>0\,\,\wrtext{and}\,\, \theta_0\in
\real\,\, \wrtext{such that}
\\ \nonumber & & \left(z_0+(0,\eps_0)e^{i(\theta_0-\eps_0,\theta_0+\eps_0)}\right)\cap
\Sigma(p)=\emptyset.
\end{eqnarray}

Our final assumption is that $p-z_0$ vanishes at $(0,0)$ precisely
to the second order, so that in a neighborhood of this point, we
have
$$
\abs{p(X)-z_0}\sim \abs{X}^2.
$$

The following is the main result of this work.

\begin{theo}
Assume that $z_0$ is at the boundary of the semiclassical
pseudospectrum of $P$, defined in {\rm (\ref{1.3.5})}, and that
$p^{-1}(z_0)\cap \real^2=(0,0)$, with $p-z_0$ vanishing there
precisely to the second order. Furthermore, we make an exterior
cone assumption {\rm (\ref{1.3.6})}. Then there exists $\alpha \in
\comp$, $\abs{\alpha}=1$, such that $\Re \alpha p''(0,0)>0$, and
the spectrum of $P$ in a sufficiently small but fixed \neigh{} of
$z_0$ in $\comp$ is given by
$$
z_k=z_0+G\left(h(k+\frac{1}{2});h\right)+{\cal
O}(h^{\infty}),\quad k\in \nat.
$$
Here $G(q;h)$ is holomorphic in $q\in {\rm neigh}(0,\comp)$, and
has an asymptotic expansion in the space of such functions, as
$h\rightarrow 0$,
$$
G(q;h)\sim \sum_{j=0}^{\infty} h^j G_{j}(q).
$$
We have $G_{0}(0)=0$ and
$$
\arg\left(\frac{\partial}{\partial q}
G_{0}(0)\right)=\arg\left(\frac{\mu}{i}\right),
$$
where $\mu$ is the unique eigenvalue of the Hamilton map of
$p''(0,0)$, for which
$$
\Re (\alpha \mu/i)>0.
$$
\end{theo}

\medskip
We illustrate Theorem 1.1 by comparing the results of a numerical
computation of small eigenvalues of the operator
$P=(hD_x)^2+cx^2+x^4$, when $c=1+3i$, with the direction of
$\arg(c)/2$, given by the quadratic part of the symbol---see
Figure 1 on the next page.

\bigskip
\noindent {\it Acknowledgement}. I am grateful to Maciej Zworski
for fruitful discussions on the subject of the present note. I
would also like to thank Caroline Lasser for helpful discussions
and advice, and Johannes Sj\"ostrand for stimulating comments. The
partial support of the National Science Foundation under grant
DMS-0304970 is gratefully acknowledged.

\begin{figure}
\begin{center}
\scalebox{.8}{\includegraphics{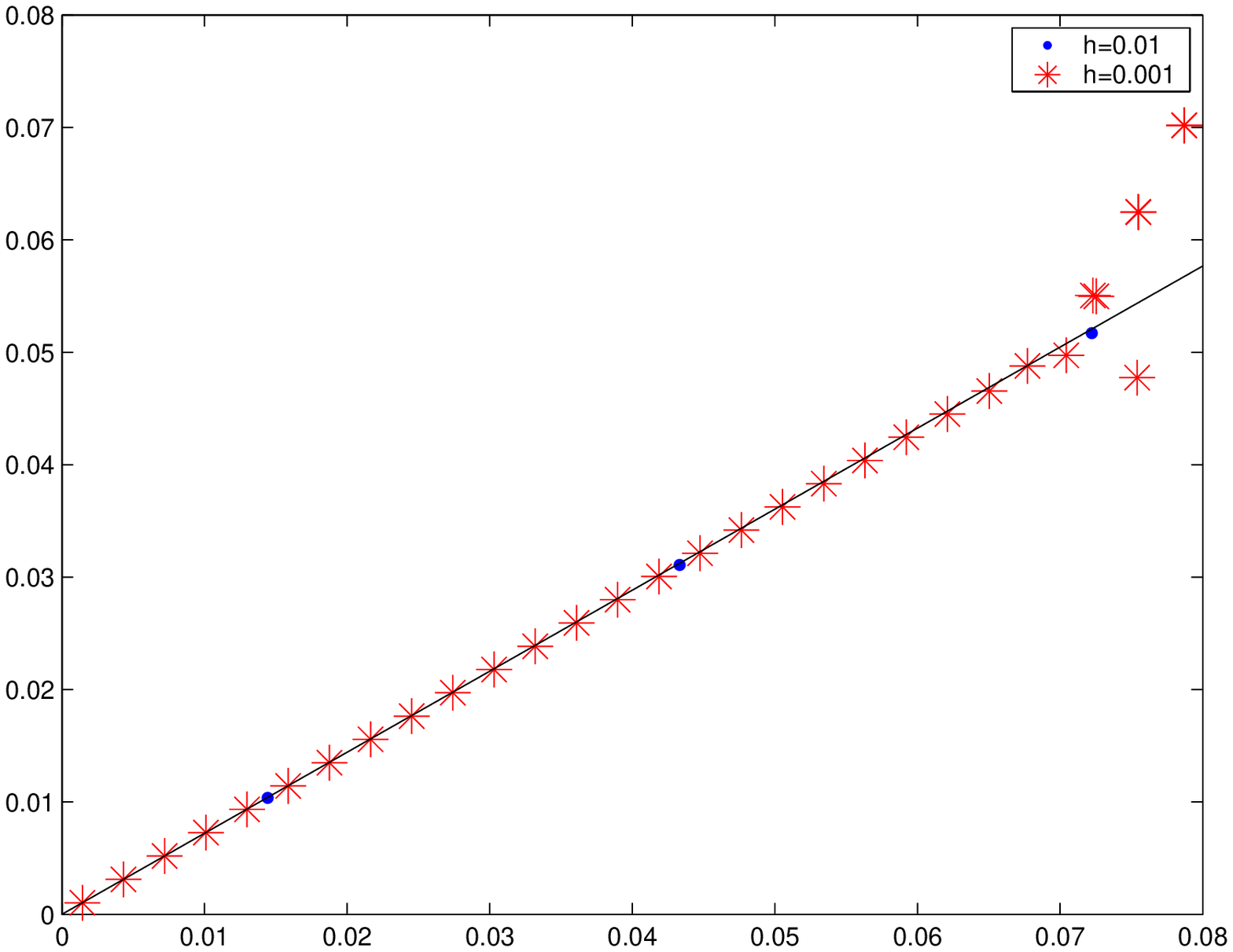}} \caption{Numerical
computation of the eigenvalues of $P=(hD_x)^2+cx^2+x^4$, when
$c=1+3i$, for two different values of $h$. When computing the
eigenvalues of $P$, following~\cite{Tref}, we discretized the
operator using the Chebyshev spectral method.} To compare the
results of the computation with Theorem 1.1 we have also plotted
the black solid line of the slope $\tan(\arg(c)/2)\approx 0.7208$.
\end{center}
\end{figure}

\section{Proof of Theorem 1.1}
\label{section1}
\setcounter{equation}{0}

Throughout this section, we shall assume, as we may, that $z_0=0$. Our starting point is the
following essentially well-known result.

\begin{prop}
Let $q(x,\xi)$ be a complex-valued quadratic form on $\real^2$
such that $q(x,\xi)=0$ precisely when $(x,\xi)=(0,0)$. Then the
range of $q$ on $\real^2$ is either all of $\comp$ or there exists
a proper closed convex cone $\subset \comp$ which contains the range
of $q$.
\end{prop}
\begin{proof}
This result is established in Lemma 3.1 of~\cite{Sj1}, and we
shall recall the proof for completeness only. In suitable linear
coordinates on $\real^2$, we may write for some $c\neq 0$,
$$
q(x,\xi)=c(x-i\xi)(x-\lambda\xi),\quad \Im \lambda\neq 0.
$$
Writing $z=x+i\xi$, we find that in the case when $\Im \lambda
<0$,
$$
q=\frac{1}{\alpha} \overline{z}(z+\gamma \overline{z}),\quad
\alpha=\frac{2i}{c(i-\lambda)},
$$
where $\gamma=(i+\lambda)(i-\lambda)^{-1}$ satisfies
$\abs{\gamma}<1$. Therefore,
$$
\Re (\alpha q)=\Re (\abs{z}^2+\gamma (\overline{z})^2)=\Re
\left(\abs{z}^2\left(1+\gamma
\frac{(\overline{z})^2}{\abs{z}^2}\right)\right)>0,
$$
when $z\neq 0$, and in this case the range of $q$ is contained in
a closed angle of opening $<\pi$. In the opposite case, $\Im \lambda>0$, we
find that with some $\beta\neq 0$,
$$
\beta q=\gamma\abs{z}^2+(\overline{z})^2,\quad \abs{\gamma}<1.
$$
It follows that the argument variation of $\beta q$ along the
circle $\abs{z}=1$ is non-vanishing, and hence the range of $q$ on
$\real^2$ is all of $\comp$.
\end{proof}

\Remark. According to Theorem 2.1.18 in~\cite{Ho1}, the range of a
general complex-valued quadratic form on $\real^2$, restricted to
$\abs{z}=1$, is an ellipse, possibly degenerated to an interval or a point.

Applying Proposition 2.1 to the quadratic form
$p_2(X)=\langle{p''(0,0)X,X}\rangle/2$, which begins the Taylor
expansion of $p$ at the origin, and using the exterior cone
condition (\ref{1.3.6}), we conclude that since $0$ is a boundary
point of the pseudospectrum of $P$, the range of $p_2$ is not all
of $\comp$, for otherwise, the argument variation of $p$ along the
positively oriented boundary of every small disc centered at $0$,
would be non-vanishing. We then know that for some $\alpha\in
\comp$, $\abs{\alpha}=1$, the real part of $\alpha p_2$ is
positive definite. After a multiplication of $P$ by $\alpha$, we
may and will assume henceforth that $\alpha=1$, so that the range
of $p_2$ on $\real^2$ is contained in a set of the form
$$
\Gamma=\{w\in \comp; \abs{\Im w}\leq C\Re w\},
$$
for some $C>0$.

By polarizaton, $p_2(X)$ gives rise to a symmetric bilinear form
on $\comp^2$, and then we write
$$
p_2(X,Y)=\sigma(X,FY), \quad X,Y\in \comp^2,
$$
where $\sigma$ is the complex symplectic form on $\comp^2$, and
$F:\comp^2 \rightarrow \comp^2$ is the Hamilton map of $p_2$---see
Section 21.5 of~\cite{Ho} and also~\cite{Ho2}, for a systematic
discussion of quadratic forms on symplectic vector spaces. Since
$\Re p_2$ is positive definite on $\real^2$, it is true that $F$
is bijective, and it is well known and easily seen that if $\mu$
is an eigenvalue of $F$ then so is $-\mu$. When $X$ is an
eigenvector of $F$ corresponding to the eigenvalue $\mu$, we write
$$
p_2(\overline{X},X)=\sigma(\overline{X},\mu X)=2i\mu \sigma(\Re
X,\Im X).
$$
Here the left hand side is equal to $p_2(\Re X,\Re X)+p_2(\Im
X,\Im X)$, and then we see that either $\mu/i \in \Gamma$ or
$-\mu/i\in \Gamma$. In what follows we shall let $\mu$ stand for
the unique eigenvalue of $F$ for which $\Re (\mu/i)>0$. Our goal
is to show that as $h\rightarrow 0$, the spectrum of $P$ near $0$
accumulates towards the origin in the direction given by the
argument of $\mu/i$.

\Remark. In Theorem 3.5 of~\cite{Sj1}, the spectrum of the Weyl
quantization of a general complex valued elliptic quadratic form
was computed. When specialized to the case at hand, the results
of~\cite{Sj1} show that the spectrum  of $p_2^w(x,hD)$ is of the
form $\{-i\mu h(2k+1),\,k=0,1,2,\ldots \}$. Our Theorem 1.1 should
therefore be viewed as an extension of the one-dimensional version
of Theorem 3.5 of~\cite{Sj1} to the case of a general analytic
$h$-pseudodifferential operator. The proof of Theorem 3.5
in~\cite{Sj1} employs complex linear canonical transformations, in
order to essentially reduce the quadratic form globally to a
complex multiple of the harmonic oscillator. In our case, we shall
proceed similarly, but now merely microlocally, working near the
critical point $(0,0)$ of the leading symbol of $P$.

\bigskip
\noindent We shall first use real canonical transformations to
simplify $p_2$. After a real linear symplectic change of
coordinates, and a conjugation of $P$ by means of a
correspon\-ding unitary metaplectic operator, we may assume that
as $(x,\xi)\rightarrow (0,0)$, \begeq \label{2.1}
p(x,\xi)=p_2(x,\xi)+{\cal
O}((x,\xi)^3)=\frac{\lambda}{2}(x^2+\xi^2)+\frac{i}{2}(ax^2+b\xi^2+2c
x\xi)+{\cal O}((x,\xi)^3), \quad \lambda>0,
\endeq
for some $a,b,c\in \real$. Here the quadratic form $\Im p_2$ can
be diagonalized by means of an orthogonal transformation, which
may be assumed to be orientation preser\-ving. Hence this
transformation is symplectic, since we are in dimension one, and
we conclude that after an additional real linear symplectic
transformation and a metaplectic conjugation, we may reduce the
quadratic part of $p$ to the form in (\ref{2.1}) with $c=0$. We
then write
$$
p(x,\xi)=\frac{(\lambda+ia)}{2}\left(x^2+d\xi^2\right)+{\cal
O}((x,\xi)^3),
$$
where $d=(\lambda+ia)^{-1}(\lambda+ib)$, and it is no restriction
to assume that $\Im d\geq 0$. After an additional real symplectic
dilation in $(x,\xi)$, the symbol $p(x,\xi)$ takes the form
\begeq
\label{2.2}
p(x,\xi)=\frac{(\lambda+ia)\abs{d}^{1/2}}{2}\left(x^2+e^{i\alpha}\xi^2\right)+{\cal
O}((x,\xi)^3),\quad \alpha=\arg d\in [0,\pi).
\endeq

Assuming that $\alpha>0$, we shall now look for a small
IR-deformation $\Lambda$ of $\real^2$, such that the restriction
of the quadratic part of $p$ to $\Lambda$ will have a constant
argument, in a neighborhood of $(0,0)$. We first argue formally,
and consider $\Lambda_0=T^*(e^{i\theta}\real)$, where $\theta\in
(0,\pi/4)$ is to be chosen. Notice that with $G(x,\xi)=\theta
x\xi$, we can write $\Lambda_0=\exp(iH_{G})(\real^2)$. If
$(x,\xi)\in \Lambda_0$ then $x=e^{i\theta}y$,
$\xi=e^{-i\theta}\eta$ for $(y,\eta)\in \real^2$, and using
$(y,\eta)$ as real symplectic coordinates on $\Lambda_0$, we get
$$
p\big|_{\Lambda_0}=p(\exp(iH_G)(y,\eta))=\frac{(\lambda+ia)\abs{d}^{1/2}}{2}\left(e^{2i\theta}y^2+
e^{i\alpha}e^{-2i\theta}\eta^2\right)+{\cal O}((y,\eta)^3).
$$
We see that if we choose $\theta=\alpha/4$, then
\begin{eqnarray*}
p(\exp(iH_G)(y,\eta)) & = &
\frac{(\lambda+ia)d^{1/2}}{2}\left(y^2+\eta^2\right)+{\cal
O}((y,\eta)^3)
\\
& = & \frac{((\lambda+ia)(\lambda+ib))^{1/2}}{2}(y^2+\eta^2)+{\cal
O}((y,\eta)^3).
\end{eqnarray*}
Here we take the square root with the positive real part. Using
the symplectic invariance of the Hamilton map $F$ of $p_2$, we then check
that here the coefficient
$$
((\lambda+ia)(\lambda+ib))^{1/2}/2
$$
is equal to $\mu/i$, where $\mu$ is the eigenvalue of $F$, for
which $\Re(\mu/i)>0$.

We shall now introduce a global IR-manifold $\Lambda$ which in a
very small complex \neigh{} of $(0,0)$ agrees with
$T^*(e^{i\theta}\real)$, and further away from that set it is
equal to $\real^2$. The manifold $\Lambda$ will also be chosen so
that the restriction of $p$ to $\Lambda\backslash\{(0,0)\}$ is
non-vanishing. We remark that the following construction of
$\Lambda$ is similar to some arguments of~\cite{KaKe}
and~\cite{Sj3} in the theory of resonances.

It will be convenient to work with the following FBI--Bargmann
transform,
\begeq
\label{2.3} Tu(x)=Ch^{-3/4} \int e^{i\varphi(x,y)/h}
u(y)\,dy,\quad x\in \comp, \quad C>0,
\endeq
where
$$
\varphi(x,y)=\frac{i}{e^{i\alpha/2}}\left(\frac{x^2}{2}+\frac{y^2}{2}-\sqrt{2}yx\right).
$$
Notice that $\Im \varphi''_{yy}>0$ and $\varphi''_{xy}\neq 0$, so
that $\varphi(x,y)$ is an admissible phase in (\ref{2.3}).
Associated to $T$, there is a complex linear canonical
transformation $\kappa_T: (y,-\varphi'_y(x,y))\mapsto
(x,\varphi_x(x,y))$, given by
$$
\kappa_T(y,\eta)=\frac{1}{\sqrt{2}}\left(y-ie^{i\alpha/2}\eta,-ie^{-i\alpha/2}y+\eta\right),\quad
(y,\eta)\in \comp^2.
$$
From the general theory~\cite{Sj2}, we know that $\kappa_T$ maps
$\real^2$ bijectively onto
$$
\Lambda_{\Phi_0}:=\left\{\left(x,\frac{2}{i}\partial_x \Phi_0(x)\right),\,x\in
\comp\right\},
$$
where $\Phi_0$ is the strictly subharmonic quadratic form given by
$$
\Phi_0(x)=\sup_{y\in {\rm {\bf R}}}-\Im \varphi(x,y).
$$
A straightforward computation shows that
\begeq
\label{2.3.5}
\Phi_0(x)=\cos(\alpha/2)\frac{(\Re
x)^2}{2}+\frac{1+\sin^2(\alpha/2)}{\cos(\alpha/2)}\frac{(\Im
x)^2}{2}+\sin(\alpha/2)\Re x\Im x.
\endeq
We remark that $\Phi_0$ is even strictly convex. Another
computation shows next that the image of $T^*(e^{i\theta}\real)$,
$\theta=\alpha/4$, under $\kappa_T$ is of the form
$\Lambda_{\Phi_1}$, where $\Phi_1(x)=\frac{1}{2}\abs{x}^2$.

When $\eps>0$, we take $\chi_{\eps}\in C^{\infty}_0(\comp)$,
$0\leq \chi\leq 1$, such that $\chi_{\eps}=1$ when
$\abs{x}<e^{1/\eps}$, $\chi_{\eps}=0$ for $\abs{x}>e^{2/\eps}$,
and such that
$$
\partial^{\alpha}\chi_{\eps}(x)={\cal
O}_{\alpha}(\eps)\abs{x}^{-\abs{\alpha}},\quad \abs{\alpha}>0.
$$
Then the smooth function
$$
\widetilde{\Phi}=\chi_{\eps}\Phi_1+(1-\chi_{\eps})\Phi_0
$$
is strictly convex for $\eps>0$ small enough, and satisfies
$$
\widetilde{\Phi}(x)=\Phi_1(x),\quad \abs{x}\leq e^{1/\eps},
$$
$$
\widetilde{\Phi}(x)=\Phi_0(x),\quad \abs{x}\geq e^{2/\eps}.
$$
In order to decrease the size of the neighborhood of $0$ where
$\widetilde{\Phi}\neq \Phi_0$, we introduce the strictly convex
function
\begeq
\label{2.3.6}
\Phi_{\eta}(x)=\left(\frac{\eta}{K_{\eps}}\right)^2
\widetilde{\Phi}\left(\frac{K_{\eps}x}{\eta}\right),\quad
K_{\eps}:=e^{2/\eps},\quad 0<\eta\ll 1,
\endeq
so that $\Phi_{\eta}(x)=\Phi_0(x)$ for $\abs{x}\geq \eta$. Notice
that taking $\eta$ sufficiently small, we may achieve that
$\Phi_{\eta}-\Phi_0$ is arbitrarily small in $C^1$--norm. We then
introduce the IR-manifold
$$
\Lambda_{\Phi_{\eta}}: \xi=\frac{2}{i}\partial_x \Phi_{\eta},
$$
and remark that along $\Lambda_{\Phi_{\eta}}$, we have
$$
\xi=\xi(x)=\chi_{\eps}\left(\frac{K_{\eps}x}{\eta}\right)\frac{2}{i}\frac{\partial
\Phi_1(x)}{\partial x}+\
\left(1-\chi_{\eps}\left(\frac{K_{\eps}x}{\eta}\right)\right)\frac{2}{i}\frac{\partial
\Phi_0(x)}{\partial x}+{\cal O}(\eps \abs{x}),
$$
uniformly in $\eta$. We claim now that with $\widetilde{p}=p\circ\kappa_T^{-1}$,
\begeq \label{2.4}
 \abs{\widetilde{p}(x,\xi(x)}\geq
\frac{1}{\widetilde{C}}\abs{x}^2,\quad \widetilde{C}>0,
\endeq
if $x\in {\rm neigh}(0,\comp)$ and $\eps>0$ is small enough.
Indeed, using (\ref{2.2}) and computing the inverse of $\kappa_T$,
we find that
$$
\widetilde{p}(x,\xi)=2((\lambda+ia)(\lambda+ib))^{1/2}ix\xi+{\cal
O}((x,\xi)^3).
$$
Using the strict convexity of the quadratic forms $\Phi_0$ and
$\Phi_1$, we get
\begin{eqnarray*}
\Re (ix\xi(x)) & = & \chi_{\eps}\abs{x}^2+(1-\chi_{\eps})\Re
(2x\partial_x \Phi_0)+{\cal O}(\eps\abs{x}^2)
\\
& = & \chi_{\eps}\abs{x}^2+(1-\chi_{\eps})\langle{x,\nabla
\Phi_0(x)}\rangle+{\cal O}(\eps \abs{x}^2)\sim \abs{x}^2,
\end{eqnarray*}
uniformly in $\eps$ and $\eta$, and the transversal ellipticity
property (\ref{2.4}) follows at once. We conclude that for every
fixed small $\eta>0$, the IR-manifold
$\Lambda:=(\kappa_T)^{-1}(\Lambda_{\Phi_{\eta}})$ is equal to
$\real^2$ away from an arbitrarily small, previously given
neighborhood of $(0,0)$ and it agrees with $T^*(e^{i\theta}\real)$
near $(0,0)$. Moreover, it is clear that the restriction of $p$ to
$\Lambda$ vanishes precisely at $(0,0)$.

Let us now recall the $h$-dependent Hilbert space
$H(\Lambda)\simeq L^2(\real)$, associated with the IR-manifold
$\Lambda$ and equipped with the norm obtained by replacing the
weight $\Phi_0$ by $\Phi=\Phi_{\eta}$ on the FBI transform side.
We also introduce the Sobolev space $H(\Lambda,m)$ which agrees
with $H(m)$ as a set, and whose norm is obtained from the
$H(m)$--norm by the same modification of the weight. Deforming the
contour in the integral representation of $P$ on the transform
side, as explained in~\cite{MeSj}, we see that \begeq \label{2.5}
P={\cal O}(1): H(\Lambda, m)\rightarrow H(\Lambda).
\endeq

We summarize the discussion above in the following proposition.

\begin{prop}
Let $P$ be as in the introduction and recall that in suitable real
symplectic coordinates, the leading symbol $p$ of $P$ has the form
{\rm (\ref{2.2})}. There exists an IR-manifold $\Lambda\subset
\comp^2$, which coincides with $\Lambda_0=T^*(e^{i\theta}\real)$,
$\theta=\alpha/4$, near $(0,0)$ and agrees with $\real^2$ outside
an arbitrarily small \neigh{} of $(0,0)$, such that $p\neq 0$
along $\Lambda\backslash\{(0,0)\}$. After applying the canonical
transformation $\kappa_T$, associated with the transform {\rm
(\ref{2.3})}, so that $\real^2$ becomes $\Lambda_{\Phi_0}$, with
$\Phi_0$ defined in {\rm (\ref{2.3.5})}, $\Lambda$ becomes
$\Lambda_{\Phi}$, where $\Phi$ is a strictly convex smooth
function defined in {\rm (\ref{2.3.6})}. For $x\in {\rm
neigh}(0,\comp)$, we have
$$
\abs{p\circ \kappa_T^{-1}\left(x,\frac{2}{i}\frac{\partial
\Phi}{\partial x}\right)}\sim \abs{x}^2.
$$
There exists a Fourier integral operator
$$
U=e^{-G(x,hD_x)/h}={\cal O}(1): H(\Lambda)\rightarrow L^2(\real),
$$
microlocally unitary near $(0,0)$, such that with $P$ denoting the
operator in {\rm (\ref{2.5})}, we have
$$
UP=\widetilde{P}U,
$$
microlocally near $(0,0)$. Here $\widetilde{P}: L^2(\real)\rightarrow
L^2(\real)$ has the leading symbol
$$
p(\exp(iH_G)(x,\xi))=\widetilde{p}(x,\xi)=
\frac{\mu}{i}\left(x^2+\xi^2\right) +{\cal O}((x,\xi)^3),\quad
(x,\xi)\rightarrow (0,0),
$$
where $\mu$ is the eigenvalue of the Hamilton map of the quadratic
part of $p$ at $(0,0)$, for which $\Im \mu>0$.
\end{prop}

\vskip 4mm In what follows we shall be interested in eigenvalues
$E$ of $P$ with $\abs{E}={\cal O}(\eps^2)$, when $0<\eps \ll 1$ is
sufficiently small but fixed. The eigenfunctions corresponding to
such eigenvalues are concentrated in a region where
$\abs{(x,\xi)}={\cal O}(\eps)$, and therefore we introduce a
change of variables $x=\eps y$. Then
$$
\frac{1}{\eps^2}P(x,hD_x;h)=\frac{1}{\eps^2}
P(\eps(y,\widetilde{h}D_y);h),\quad \widetilde{h}=\frac{h}{\eps^2}.
$$
When viewed as an $\widetilde{h}$-pseudodifferential operator,
$\eps^{-2}P$ has a corresponding new symbol
$$
\frac{1}{\eps^2} P(\eps(y,\eta))\sim \frac{1}{\eps^2}
p(\eps(y,\eta))+\widetilde{h}p_1(\eps(y,\eta))+\ldots,
$$
to be considered in a region where $\abs{(y,\eta)}={\cal O}(1)$.

This scaling reduction together with Proposition 2.2 allows us to
reduce the further analysis to an analytic $\widetilde{h}$-pseudodifferential
operator $P_{\eps}$, microlocally defined near $(0,0)\in \real^2$,
whose (not necessarily real) leading symbol has the form
$$
p_{\eps}(x,\xi):=p_{0,\eps}(x,\xi)=\frac{x^2+\xi^2}{2}+\eps r_{\eps}(x,\xi),\quad
r_{\eps}(x,\xi)={\cal O}((x,\xi)^3).
$$
The full symbol is
$$
P_{\eps}(x,\xi;\widetilde{h})\sim \sum_{j=0}^{\infty}
\widetilde{h}^j p_{j,\eps}(x,\xi),
$$
with $p_{j,\eps}(x,\xi)$ holomorphic in a fixed complex \neigh{} of
$(x,\xi)=(0,0)$.

From Section 5 of~\cite{HiSj1} we recall that for every small
enough $\eps$, there exists a holomorphic canonical transformation
$$
\kappa_{\eps}: {\rm neigh}((0,0),\comp^2)\rightarrow {\rm
neigh}((0,0),\comp^2),
$$
such that $\kappa_{\eps}(0,0)=(0,0)$, and $\Im
\kappa_{\eps}(x,\xi)={\cal O}(\eps)$ when $x,\xi$ are real, and
such that
$$
p_{\eps}\circ \kappa_{\eps}(x,\xi)=g_{\eps}\left(\frac{x^2+\xi^2}{2}\right).
$$
Here $g_{\eps}(E)$ is an analytic function of $E$, with
$g_{\eps}(0)=0$, and such that $\Im g_{\eps}(E)={\cal O}(\eps)$
for real $E$, and $\Re g_{\eps}'(0)>0$. We may also remark that we
have in fact $g_{\eps}(E)=E+{\cal O}(\eps E)$.

Associated to the canonical transformation $\kappa_{\eps}$, we
introduce a global IR-manifold $\Lambda_{\eps}\subset \comp^2$,
which is $\eps$--close to $\real^2$, equals
$\kappa_{\eps}(\real^2)$ in a complex \neigh{} of $(0,0)$, and
agrees with the real phase space further away from this set.
Precisely as in~\cite{HiSj1}, we then implement $\kappa_{\eps}$ by
means of an elliptic $\widetilde{h}$--Fourier integral operator
$U_{\eps}$,
$$
U_{\eps}={\cal O}(1): L^2(\real)\rightarrow H(\Lambda_{\eps}).
$$
Then the action of $P_{\eps}$ on $H(\Lambda_{\eps})$ is, microlocally
near $(0,0)$, unitarily equivalent to the
$\widetilde{h}$-pseudodifferential operator
$$
\widetilde{P}_{\eps}=U^{-1}_{\eps} P_{\eps} U_{\eps}:
L^2(\real)\rightarrow L^2(\real),
$$
whose complete Weyl symbol has the form
$$
\widetilde{P}_{\eps}(x,\xi;\widetilde{h})\sim \sum_{j=0}^{\infty}
\widetilde{h}^j \widetilde{p}_{j,\eps}(x,\xi),
$$
with
$$
\widetilde{p}_{0,\eps}(x,\xi)=g_{\eps}\left(\frac{x^2+\xi^2}{2}\right).
$$

Continuing to follow~\cite{HiSj1} (where also further references
are given), let us now recall how to simplify the lower order
terms in an operator whose principal symbol is modelled on the
one-dimensional harmonic oscillator. In doing so, we consider a
formal $h$-pseudodifferential operator $Q(x,hD_x;h)$, defined
microlocally near $(0,0)\in \real^2$, with symbol
$$
Q(x,\xi;h)\sim q_0(x,\xi)+h q_1(x,\xi)+\ldots,
$$
where $q_j$ are holomorphic in a fixed complex \neigh{} of
$(x,\xi)=(0,0)$, and
$$
q_0(x,\xi)=g_0((x^2+\xi^2)/2),
$$
with $g_0$ holomorphic near $0$, $g_0(0)=0$, $g'_0(0)\neq 0$. In
Section 5 of~\cite{HiSj1} we recalled that by means of an
averaging procedure we can construct
$$
A(x,\xi;h)\sim a_0(x,\xi)+h a_1(x,\xi)+\ldots,
$$
with all $a_j$ holomorphic in a $j$-independent \neigh{} of
$(0,0)$, such that
\begeq \label{2.6}
 e^{iA(x,hD_x;h)} Q(x,hD_x;h)
e^{-iA(x,hD_x;h)}=g\left(\frac{x^2+(hD_x)^2}{2};h\right).
\endeq
Here
$$
g(E;h)\sim \sum_{j=0}^{\infty} g_j(E) h^{j},
$$
with $g_0$ as above, and $g_j(E)$ holomorphic in a fixed \neigh{}
of $E=0$. As in~\cite{HiSj1}, when deriving (\ref{2.6}) we make use of the fact
that if $f$ is a smooth (holomorphic) function on $\real^2$ which
depends on $(x^2+\xi^2)/2$ only, then $f^w(x,hD_x)=\widetilde{f}\left(\frac{1}{2}(x^2+(hD_x)^2)\right)$,
for some $\widetilde{f}$ with $\widetilde{f}=f+{\cal
O}(h)$---see~\cite{Ho2} for this essentially well-known result.

Applying this discussion to $\widetilde{P}_{\eps}$, we conclude
that there exists
$$
A_{\eps}(x,\xi;\widetilde{h})\sim
A_{0,\eps}(x,\xi)+\widetilde{h}A_{1,\eps}(x,\xi)+\ldots,
$$
with each $A_{j,\eps}$ holomorphic in a fixed \neigh{} of $(0,0)$ in
$\comp^2$, such that after a conjugation by
$e^{iA_{\eps}(x,\widetilde{h}D_x;\widetilde{h})}$, the operator
$\widetilde{P}_{\eps}$ becomes
$$
\widehat{P}_{\eps}=G_{\eps}\left(\frac{x^2+(\widetilde{h}D_x)^2}{2};\widetilde{h}\right),
$$
where
$$
G_{\eps}(q;\widetilde{h})\sim \sum_{j=0}^{\infty} \widetilde{h}^j
G_{j,\eps}(q),
$$
with all $G_{j,\eps}$ holomorphic in $q$ in a $j$-independent \neigh{} of
$0\in \comp$. Moreover, $G_{0,\eps}=g_{\eps}$ above.

Combining Proposition 2.2 together with the reductions above, we
arrive at the following result.

\begin{prop}
Let $P$ and $p$ be as in the introduction and recall that we
assume that the real part of the Hessian of $p$ at $(x,\xi)=(0,0)$ is
positive definite. Let $\eps>0$ be small enough but fixed, and let
us view the operator
$$
\frac{1}{\eps^2}P(x,hD_x)=\frac{1}{\eps^2}
P\left(\eps(y,\widetilde{h} D_y)\right),\quad x=\eps y,\quad
\widetilde{h}=\frac{h}{\eps^2},
$$
as an $\widetilde{h}$-pseudodifferential operator, with the
leading symbol $p_{\eps}(y,\eta)=\frac{1}{\eps^2}
p(\eps(y,\eta))$. There exists an IR-manifold
$\widehat{\Lambda}\subset \comp^2$, with $(0,0)\in
\widehat{\Lambda}$, which coincides with $\real^2$ away from a
\neigh{} of $(0,0)$ and is close to $\real^2$ everywhere, such
that the restriction of $p_{\eps}$ to $\widehat{\Lambda}$ is
elliptic outside an arbitrarily small \neigh{} of $(0,0)$.
Furthermore, there exists an elliptic Fourier integral operator
$$
U: H(\widehat{\Lambda})\rightarrow L^2(\real),
$$
such that, microlocally near $(0,0)\in \widehat{\Lambda}$,
$$
U\left(\frac{1}{\eps^2} P\right)=\widehat{P}_{\eps} U.
$$
Here
$\widehat{P}_{\eps}=\widehat{P}_{\eps}\left((1/2)((\widetilde{h}D_x)^2+x^2);\widetilde{h}\right)$
has the complete Weyl symbol,
$$
\widehat{P}_{\eps}(x,\xi;\widetilde{h})\sim \sum_{j=0}^{\infty}
\widetilde{h}^j G_{j,\eps}\left(\frac{x^2+\xi^2}{2}\right),
$$
with the functions $G_{j,\eps}$ holomorphic in a common \neigh{}
of $0\in \comp$. We have $G_{0,\eps}(0)=0$ and
$$
\arg\left(G'_{0,\eps}(0)\right)=\arg\left(\frac{\mu}{i}\right)+{\cal
O}(\eps),
$$
where $\mu$ is the unique eigenvalue of the Hamilton map of the
Hessian of $p$ at $(0,0)$, for which $\Im \mu>0$.
\end{prop}

\vskip 2mm \noindent With Proposition 2.3 available, we are now in
a position to solve a suitable Grushin problem in the globally
defined Hilbert space $H(\widehat{\Lambda})$, and it is clear that
the setup of the relevant Grushin problem is exactly the same as
in~\cite{HiSj1}. (See also~\cite{SjZw} for a recent general
presentation of Grushin problem techniques in semiclassical
analysis and spectral theory.) The analysis of the Grushin problem
in question is even simplified when compared with~\cite{HiSj1}, as
we can use ellipticity along $\widehat{\Lambda}$, when away from a
\neigh{} of $(0,0)$, while near $(0,0)$, the operator is reduced
microlocally to a function of the harmonic oscillator. In
conclusion, we obtain the following result.

\begin{prop} Let $\eps>0$ be small enough but independent of $h$. The
eigenvalues of $P$ in an open disc centered at $0$ with radius
${\cal O}(\eps^2)$ are of the form 
$$
\sim \sum_{j=0}^{\infty} h^j \eps^{2-2j}
G_{j,\eps}\left(\frac{h}{\eps^2}\left(k+\frac{1}{2}\right)\right),\quad
k\in \nat,
$$
with $G_{j,\eps}(q)$ holomorphic in $q\in {\rm neigh}(0,\comp)$.
We have
$$
G_{0,\eps}(q)=\frac{\mu}{i}q(1+{\cal O}(\eps)).
$$
\end{prop}

Combining Proposition 2.4 together with a scaling argument as
in~\cite{MeSj},~\cite{Sj3}, and~\cite{HiSj1}, we obtain Theorem
1.1.

\end{document}